\documentclass[12pt]{article}

\usepackage{amssymb}

\begin{document}
\newtheorem{theorem}{Theorem}[section]
\newtheorem{lemma}[theorem]{Lemma}
\newtheorem{corollary}[theorem]{Corollary}
\newtheorem{definition}[theorem]{Definition}
\newtheorem{proposition}[theorem]{Proposition}
\newtheorem{defprop}[theorem]{Definition-Proposition}
\newtheorem{example}[theorem]{Example}
\newtheorem{remark}[theorem]{Remark}
\newcommand{\Proof}{\noindent{\bf Proof:} }
\catcode`\@=11
\@addtoreset{equation}{section}
\catcode`\@=12
\renewcommand{\theequation}{\arabic{section}.\arabic{equation}}
\def\lu{\rightharpoonup}
\def\sqr#1#2{{\vcenter{\vbox{\hrule height.#2pt\hbox{\vrule width.#2pt
  height#1pt \kern#1pt \vrule width.#2pt}\hrule height.#2pt}}}}
\def\square{\mathchoice\sqr64\sqr64\sqr{2.1}3\sqr{1.5}3}
\def\argh{\rightharpoonup}
\def\m#1{m_{(#1)}}
\def\g#1{g_{(#1)}}
\def\h#1{h_{(#1)}}
\def\k#1{k_{(#1)}}
\def\M{{\cal M}}
\def\Z{{\mathbb Z}}
\def\F{{\cal F}}
\def\mm#1{m_{(#1*)}}
\def\a#1{a_{(#1)}}
\def\aa#1{a_{(#1*)}}
\def\b#1{b_{(#1)}}
\def\bb#1{b_{(#1*)}}
\def\l#1{l_{(#1)}}
\def\n#1{n_{(#1)}}
\def\nn#1{n_{(#1*)}}
\def\totimes{\tilde{\otimes}}
\def\c#1{c_{(#1)}}
\def\cc#1{c_{(#1*)}}
\def\id{{\rm id}}
\def\ydh{{\cal YD}(H)}
\def\yddh{{\cal YD}(H^*)}
\def\lefth{_{H}{\cal M}}
\def\lefthd{_{H^*}{\cal M}}
\def\righth{{\cal M}^H}
\def\righthd{{\cal M}^{H^*}}
\def\C{{\cal C}}
\def\D{{\cal D}}
\def\R{{\cal R}}
\def\U{{\cal U}}
\def\rsb{r_{s,\beta}}
\def\rsg{r_{s,\gamma}}
\def\f{{\bf k}}
\def\km{\f^{*}}
\def\llazy{Z^l_L(H)}
\def\lazy{Z_L(H)}
\def\pf{{\bf Proof:} }
\def\cl{C_{VL}(H)}
\def\cll{C_L(H)}
\def\ch{C(H)}
\def\bl{B_L(H)}
\def\vl{B_{VL}(H)}
\def\coh{H^2_L(H)}
\def\aut{Aut_{Hopf}(H)}
\def\lazyn{Z_L^n(H)}

\title{When is a cleft extension $H$-Azumaya?}
\author{Giovanna Carnovale}

\maketitle
\begin{abstract}{We show which 
$H^{op}$-cleft extensions $\f\# _\sigma H^{op}$ of $\f$ for a dual
    quasitriangular Hopf algebra $(H, r)$ are
$H$-Azumaya. The result is given in terms of bijectivity of a
 map $\theta_\sigma\colon H\to H^*$ defined in 
terms of $r$ and the $2$-cocycle $\sigma$, generalizing a well-known result for the
commutative and cocommutative case.
We illustrate the Theorem with an explicit computation for the Hopf algebras
of type $E(n)$.}
\end{abstract}



\section*{Introduction}

In the last fifty years several generalizations of the classical
Brauer group of a field $\f$ were introduced. A key r{o}le in
algebraic geometry is played by 
the Brauer group of suitable sheaves of algebras over a scheme $X$
(\cite{gro}). In representation theory and in ring theory 
generalizations are obtained either relaxing the requirements on
$\f$ and/or on the algebras, 
or by adding extra structures. In the first case one can, for
instance,  
allow $\f$ to be a commutative ring
 (\cite{aus}) or allow algebras without unit
(\cite{taylor}). In the second case one can
work in categories whose objects carry more structure, such as
a grading 
by an abelian group (\cite{wall},
\cite{LO1}), a commutative and cocommutative Hopf algebra action
(\cite{Long}), etcetera. Most of the known generalizations
are examples of the Brauer group of a symmetric monoidal category (\cite{par}).
The Brauer group of a general Hopf algebra (\cite{CVZ}) and later, 
the Brauer group of a
braided monoidal category (\cite{VZ}) seem to be the highest level of
generality 
for a Brauer group in this setting 
that has been reached so far. The most popular examples of
a braided monoidal category are the category of left modules over a
quasitriangular Hopf algebra and, dually, the category of right
comodules of a dual quasitriangular Hopf algebra. In these cases the
elements of the corresponding Brauer group are equivalence classes of
particular 
module and comodule algebras, respectively. The explicit
computation of the corresponding Brauer group is in general far from being
trivial. When studying these groups 
it is quite natural to wonder whether
well-known modules and comodules algebras 
can be seen as representatives of elements of these
Brauer groups, that is, whether they are Azumaya algebras in the
corresponding category. 
A very well-known family of comodule algebras is
provided by cleft extensions of the base field (\cite{doitake}). 

Cleft extensions are a natural generalization of Galois extensions. In
the finite-dimensional case the notion of cleft extension coincides
with the notion of a Hopf-Galois extension of a Hopf algebra. 
In the general case the notion of cleft extension is stronger
than the notion of a Hopf-Galois extension because the existence of a
normal basis is required. The theory of cleft
extensions has reached a satisfactory description in terms of crossed
systems (see \cite{doitake}, \cite{bcm}). Cleft extensions of the base
field correspond to equivalence classes of $2$-cocycles and they are
isomorphic, as algebras, to twists of $H$ by the corresponding
cocycle. 
A Theorem in \cite[\S 12.4]{caen} expresses
explicitely when the twist of a commutative, cocommutative Hopf
algebra by a Sweedler 
cocycle is $H$-Azumaya. The aim of the present paper is the
generalization of this theorem to the dual quasitriangular case 
and this is reached in Theorem \ref{analogo}. In analogy to the
cocommutative case the  
result states that such a twist is 
Azumaya in the category if and only if the linear map
$\theta_\sigma\colon H\longrightarrow H^*$ defined by
$h\mapsto(\sigma\tau*r*\sigma^{-1})(-\otimes h)$ is invertible. After
the proof of Theorem 
\ref{analogo} we discuss a dual version of the result obtaining
Corollary \ref{duanalogo}. The main result is illustrated
with an example: for every universal $r$-form $r$ of the Hopf algebas
of type $E(n)$ introduced in \cite{bdg} 
we characterize in Proposition \ref{wk} which twists of $E(n)$ are
Azumaya with respect to  
$r$. Finally, in Proposition \ref{deco}, we relate this description to
the computation of the Brauer 
group $BM(\f, E(n), R_0)$ of the category of modules of $E(n)$ with
respect to the $R$-matrix $R_0$ obtained in \cite{GioJuan3}. 

\section{Cleft extensions and $H$-Azumaya algebras}

Unless otherwise stated $H$ will denote a finite-dimensional Hopf
algebra over a field $\f$, with coproduct $\Delta$ and antipode $S$. 
All modules, comodules and algebras will be
assumed to be over $\f$, as well as unadorned tensor products. The
standard flip map $V\otimes W\to W\otimes V$ will be denoted by
$\tau$. For coproduct and right comodule structures we shall use the
notations
$\sum \n1\otimes\n2$ and $\sum \n0\otimes\n1$, respectively.  
 
The Brauer group of a braided monoidal category was defined in
\cite{VZ}. It is well-known that  if  
$H$ is a dual quasitriangular Hopf algebra with
universal $r$-form $r$ the category $\M^H$ of finite-dimensional right $H$-comodules 
is braided monoidal with $H$-comodule structure on the tensor product
given by 
$$\rho(m\otimes n)=\sum \m0\otimes\n0\otimes\n1\m1$$
and 
braiding $\psi$ given by
$\psi(m\otimes n)=\sum \n0\otimes\m0 r(\n1\otimes\m1)$ for every pair
of comodules $M$ and $N$ and every $m\in M$ and $n\in N$. In this
particular setting, an algebra in the category $\M^H$ is an $H^{op}$-comodule
algebra $A$.
In particular, if $P$ is a $H$-comodule
then $End(P)$ with the usual composition of endomorphisms and with 
comodule structure given by 
\begin{equation}\label{endo}
\rho(f)(a)=\sum f(\a0)_{(0)}\otimes S^{-1}(\a1)f(\a0)_{(1)} 
\end{equation} for every $f\in End(P)$
and every $a\in P$ is an algebra in $\M^H$. Similarly,
$End(P)^{op}$ is an $H^{op}$-comodule algebra with respect to the structure:  
\begin{equation}\label{endop}
\rho(f)(a)=\sum f(\a0)_{(0)}\otimes f(\a0)_{(1)}S(\a1).
\end{equation}
The opposite algebra $\overline A$ of an algebra $A$ in the category 
is equal to $A$ as a
$H$-comodule but its product is given by $a\circ b=\sum
\b0\a0r(\b1\otimes\a1)$. It is again an $H^{op}$-comodule algebra. 
Given two algebras $A$ and $B$ in the category we endow the comodule
$A\otimes B$ with the product
$$(a\#b)(c\#d)=\sum (a\c0\#\b0 d)\; r(\c1\otimes \b1)$$ for every
$a,\,c\in A$ and every $b,\,d\in B$. The resulting algebra is a
$H^{op}$-comodule algebra, denoted by $A\#B$.   
An algebra $A$ in $\M^H$ is called Azumaya, or $(H, r)$-Azumaya if the
$H^{op}$-comodule algebra maps
$$
\begin{array}{rl}
F\colon A\#\overline{A}&\longrightarrow End(A)\\
F(a\#\bar{b})(c)&=\sum a\c0\b0\;r(\c1\otimes\b1)
\end{array}
$$
and
$$
\begin{array}{rl}
G\colon \overline{A}\# A &\longrightarrow End(A)^{op}\\
G(\bar{a}\#b)(c)&=\sum r(\a1\otimes\c1)\a0\c0 b
\end{array}
$$
are isomorphisms. The opposite algebra of an $(H,r)$-Azumaya algebra and the
product $\#$ of two $(H,r)$-Azumaya algebras are again $(H,r)$-Azumaya algebras. 
The elements of the Brauer group $BC(\f,H,r)$ of the category
$\M^H$ are the 
equivalence classes of $(H,r)$-Azumaya algebras with respect to the
equivalence relation:
$A\sim B$ if $A\#End(P)\cong B\# End(Q)$ for some $H$-comodules $P$
and $Q$. The product $\#$ induces on $BC(\f, H, r)$ a group structure
and the inverse of a class represented by an algebra $A$ is the class
represented by the opposite algebra $\overline A$.

\smallskip
 
Dually, if $H$ is a quasitriangular Hopf algebra with
$R$-matrix $R=\sum R^{(1)}\otimes R^{(2)}$, the category $_H\M$ of
finite-dimensional left $H$-modules is monoidal, with usual $H$-module structure on the
tensor product of two modules. An algebra in the category $_H\M$ is
just an $H$-module
algebra. If $P$ is a $H$-module
then $End(P)$ with the usual composition of endomorphisms and with 
module structure given by 
\begin{equation}\label{endo2}
(h \cdot f)(m)=\sum h_{(1)} \cdot f(S(h_{(2)}) \cdot m)
\end{equation} 
for every $h\in H$, every $f\in End(P)$
and every $m\in P$ is an algebra in $_H\M$. Similarly,
$End(P)^{op}$ is also an algebra in $_H\M$ if we endow it with the
module structure:
\begin{equation}\label{endop2}
(h \cdot f)(m)=\sum h_{(2)} \cdot f(S^{-1}(h_{(1)}) \cdot m).
\end{equation}
The opposite algebra $\overline A$ of an algebra $A$ in $_H\M$ 
is equal to $A$ as an
$H$-module and its product is given by $a\circ b=\sum
(R^{(2)}.b)(R^{(1)}.a)$. It is again an $H$-module algebra. 
The tensor
product of two algebras $A$ and $B$ in the category $_H\M$ is
$A\# B\cong A\otimes B$ as modules with the
multiplication:
$$(a\#b)(c\#d)=\sum a(R^{(2)}.c)\# (R^{(1)}).b d $$ for every
$a,\,c\in A$ and every $b,\,d\in B$. An algebra $A$ in $_H\M$ is called Azumaya, or $(H, R)$-Azumaya if the
$H$-module algebra maps 
$$
\begin{array}{rl}
F'\colon A\#\overline{A}&\longrightarrow End(A)\\
F'(a\#\bar{b})(c)&=\sum a (R^{(2)}.c)(R^{(1)}.b)
\end{array}
$$
and
$$
\begin{array}{rl}
G'\colon \overline{A}\# A &\longrightarrow End(A)^{op}\\
G'(\bar{a}\#b)(c)&=\sum (R^{(2)}.a)(R^{(1)}.c) b
\end{array}
$$
are isomorphisms.
The elements of the Brauer group $BM(\f,H,R)$ of the category $_H\M$ are the
equivalence classes of $(H,R)$-Azumaya algebras with respect to the
equivalence relation:
$A\sim B$ if $A\#End(P)\cong B\# End(Q)$ for some $H$-modules $P$
and $Q$. The product in $BM(\f, H, R)$ is induced by the product $\#$,
with inverse represented by the opposite algebra. 

Computations of $BM(\f,H, R)$ have
 been carried out only in a few cases, namely:
for Sweedler's Hopf algebra $H_4$ with respect to the 
$R$-matrix $R_0$ in \cite{yinhuo} and for the remaining 
$R$-matrices in \cite{gio}; for the Hopf algebras of
type $H_\nu$ and all $R$-matrices in
\cite{GioJuan1}, for the group algebra of the dihedral group in
 \cite{GioJuan2}
 and for the Hopf algebras of type 
$E(n)$ and all triangular $R$-matrices in \cite{GioJuan3}.    
A key r{o}le in these computations was played by $H$-cleft extensions of
 the base field $\f$. 

\smallskip

An $H$-cleft extension $B$ of $\f$ is a
right $H$-comodule algebra such that $B^{co(H)}=\f$ 
and such that there exists a convolution invertible
 map $\gamma\colon H\to B$ (cfr. \cite{doitake}). It is well-known that
 cleft extensions of $\f$ are parametrized by $2$-cocycles, i.e.,
 convolution invertible elements $\sigma$ of $(H\otimes H)^*$
 satisfying the relations:
$$\sigma(h\otimes 1)=\sigma(1\otimes h)=\varepsilon(h)$$
$$
\sum\sigma(\k1\otimes\m1)\sigma(h\otimes\k2\m2)
=\sum\sigma(\h1\otimes\k1)\sigma(\h2\k2\otimes m)
$$
for every $h,\,k,\,m\in H$.
The cleft extension corresponding to $\sigma$ is isomorphic to the
crossed product  $_\sigma H=\f\#_\sigma H$ that is: the comodule
algebra coinciding with $H$ as a comodule and with product given by $h\cdot
k=\sum\sigma(\h1\otimes\k1)\h2\k2$. 

Two cocycles $\sigma$ and $\omega$ 
are called cohomologous if there
exists a convolution invertible element $\theta$ in $H^*$ for which
$$\sigma(h\otimes
k)=\omega^\theta(h\otimes k)=\sum\theta(\h1)\theta(\k1)
\omega(\h2\otimes\k2)\theta^{-1}(\h3\k3).$$ Two cleft extensions are
equivalent if and only if they correspond to   
cohomologous cocycles.

The cleft extension $_\sigma H$ is also a left comodule algebra for
Doi's twisted Hopf algebra $_\sigma H_{\sigma^{-1}}$. The latter is
obtained with the procedure dual to Drinfeld's twist and it is equal to
$H$ as a coalgebra but with product:
$$h\cdot_\sigma k=\sum\sigma(\h1\otimes\k1)\h2\k2\sigma^{-1}(\h3\otimes\k3).$$
It is well-known that if $H$ is dual quasitriangular with universal
$r$-form $r$ then $(\sigma\tau)*r*\sigma^{-1}$ is a universal $r$-form
for $_\sigma H_{\sigma^{-1}}$. By the dual version of \cite[Proposition 2.3.5]{mabook}, cohomologous
cocycles yield isomorphic twisted Hopf algebras and if $H$ is dual
quasitriangular, they yield isomorphic dual quasitriangular Hopf algebras.

\section{The Main result}

Given a Hopf algebra $H$, its opposite algebra $H^{op}$ with its coproduct
$\Delta$ is a left and right $H^{op}$-comodule
algebra. If $\sigma$ is a 
left $2$-cocycle for $H$ then $\sigma\tau$ is a left $2$-cocycle for
$H^{op}$ and $A_\sigma:=\;_{\sigma\tau}H^{op}$ is
again a right $H^{op}$-comodule algebra, with product:
%
$$
h\cdot k=\sum \sigma(\k1\otimes \h1)\k2\h2
$$
for every $h$ and $k\in A_\sigma$, so that
$$hk=\sum \sigma^{-1}(\h1\otimes \k1)\k2\cdot\h2.$$

If $H$ is dual quasitriangular with universal $r$-form $r$, the
product $\bullet$ in $\overline{A}_\sigma$, the opposite algebra with
respect to $r$, is given by 
$$
\begin{array}{rl}
a\bullet b&=\sum \b1\cdot\a1 r(\b2\otimes \a2) \\
&=\sum((\sigma\tau)*r)(\b1\otimes\a1)\b2\a2.
\end{array}
$$
One may wonder
when $A_\sigma$ is an $(H, r)$-Azumaya algebra. In this particular
case, the maps $F$ and $G$ in Section 1  are: 
$$
\begin{array}{rl}
F(h\# k)(l)&=\sum h\cdot (\l1\cdot \k1)r(\l2\otimes\k2)\\
&=\sum h\cdot \sigma(\k1\otimes\l1)\k2\l2 r(\l3\otimes\k3)\\
&=\sum h\cdot \sigma(\k1\otimes\l1) r(\l2\otimes\k2)\l3\k3\\ 
&=\sum h\cdot (\sigma\tau *r)(\l1\otimes\k1)\sigma^{-1}(\l2\otimes\k2)\k3\cdot\l3\\
&=\sum(\sigma\tau * r *\sigma^{-1})(\l1\otimes\k1) h\cdot \k2\cdot\l2 
\end{array}
$$
 and
$$
\begin{array}{rl}
G(h\#k)(l)&=\sum r(\h2\otimes\l2)\h1\cdot\l1\cdot k\\
&=\sum \sigma(\l1\otimes\h1)(\l2\h2) r(\h3\otimes\l3)\cdot k\\
&=\sum \sigma(\l1\otimes\h1)r(\h2\otimes\l2)(\h3\l3)\cdot k\\
&=\sum (\sigma\tau * r*\sigma^{-1})(\h1\otimes\l1)\l2\cdot\h2\cdot k.
\end{array}
$$
The bijectivity of these maps is strictly related to the behaviour of
the universal $r$-form $r_\sigma=(\sigma\tau)*r*\sigma^{-1}$ in the
twisted Hopf algebra
$_\sigma H_{\sigma^{-1}}$. It is well-known that if $r$ is a
universal $r$-form for $H$, the map
$$
\begin{array}{rl}
\theta_r\colon H^{op}&\longrightarrow H^*\\
h&\mapsto r(-\otimes h)
\end{array}
$$
is a Hopf algebra homomorphism. In particular we will relate the
bijectivity of the map
$$
\begin{array}{rl}
\theta_{\sigma}=\theta_{r_\sigma}\colon
(\,_\sigma H_{\sigma^{-1}})^{op}&\longrightarrow (\,_\sigma
H_{\sigma^{-1}})^{*}\\
h&\mapsto r_\sigma(-\otimes h)
\end{array}
$$  
to the bijectivity of
$F$ and $G$. We shall follow the lines of the proof of \cite[Theorem 12.4.5]{caen}. In terms of $\theta_\sigma$ we have:
$$F(h\#k)(l)=h\cdot\left(\sum\langle  \theta_\sigma(\k1),\,\l1\rangle
\k2\cdot\l2\right)$$ 
$$G(h\#k)(l)=\left(\sum\langle \theta_\sigma(\l1),\,\h1\rangle
\l2\cdot\h2\right)\cdot k.$$ 

%

\smallskip

We recall that for a finite-dimensional Hopf algebra $H$ the space of
left integrals $\int^l_{H^*}$ for $H^*$ is one-dimensional $\f\zeta$,
say. As a consequence of the Fundamental Theorem for Hopf modules 
there is a $\f$-linear isomorphism  
$$
\begin{array}{rl}
V\colon \int_{H^*}^l\otimes H&\longrightarrow H^*\\
V(\xi\otimes h)(k)&=\xi(kS(h)).
\end{array}
$$

It is well-known that if we put $v(h):=V(\zeta\otimes h)$ the
following formula holds: 

\begin{equation}\label{pinco}
\sum\langle v(h),\,\k2\rangle \k1=\sum\langle v(\h1),\,k\rangle \h2
\end{equation}
for every $h,\,k\in H$. 

Let us denote by $w(h):=(S^{-1})^*(v(h))\in H^*$ for every $h\in H$. Then one has:

\begin{equation}\label{pallino}
\sum\langle w(h),\,\k1\rangle S^{-1}(\k2)=\sum \langle w(\h1),\,k\rangle \h2.
\end{equation}
Applying the antipode $S$ on both sides we get:
\begin{equation}\label{vu}
\sum\langle w(h),\,\k1\rangle \,\k2=\sum\langle w(\h1),\,k\rangle \,S(\h2)
\end{equation}
which is the counterpart of (\ref{pinco}) for $H^{op,cop}$.

We introduce the following maps $S_i\colon
A_{\sigma}\longrightarrow A_{\sigma}$ for $i=1,2$: 
%
%
$$S_1(h)=\sum\sigma^{-1}(S(\h2)\otimes\h3)S(\h1);$$
%
$$S_2(h)=\sum\sigma^{-1}(\h3\otimes S^{-1}(\h2))S^{-1}(\h1).$$
%
%
%
A straightforward computation yields:
\begin{equation}\label{vare}
\sum \h2\cdot S_1(\h1)=\varepsilon(h)=\sum S_2(\h1)\cdot \h2
\end{equation}
for every $h\in A_{\sigma}$.
Besides, by the left cocycle condition we have:
\begin{equation}\label{giro}\sigma(k\otimes
lm)=\sum\sigma^{-1}(\l1\otimes\m1)
\sigma(\k1\otimes\l2)\sigma(\k2\l3\otimes\m2) 
\end{equation}
and
\begin{equation}\label{giro2}\sigma(kl\otimes
m)=\sum\sigma^{-1}(\k1\otimes\l1)\sigma(\l2\otimes\m1)\sigma(\k2\otimes\l3\m2). 
\end{equation} 
Applying (\ref{giro}) to $k=\h1$, $l=S(\h2)$ and $m=\h3$ and adding
all terms we get:
$$\varepsilon(h)=\sum\sigma^{-1}(S(\h3)\otimes\h4)
\sigma(\h1\otimes S(\h2)).$$
This formula was already observed, in greater generality, by Blattner
and Montgomery, see \cite[Proposition 7.2.7]{M}. It implies that
\begin{equation}\label{vare2}
\sum S_1(\h2)\cdot \h1=\varepsilon(h)
\end{equation}
for every $h\in A_\sigma$. Applying (\ref{giro2}) to $k=\h3$,
$l=S^{-1}(\h2)$ and $m=\h1$ and adding all terms we have:
$$\varepsilon(h)=\sum\sigma^{-1}
(\h4\otimes S^{-1}(\h3))\sigma(S^{-1}(\h2)\otimes\h1)$$
and this implies that 
\begin{equation}\label{vare3}
\sum \h1\cdot S_2(\h2)=\varepsilon(h).
\end{equation}

Noe we are ready to state the main result of this section.
\begin{theorem}\label{analogo} Let $H$ be finite-dimensional  dual
quasitriangular Hopf algebra $H$ with 
universal $r$-form $r$. Let $\sigma$
be a left $2$-cocycle for $H$. Then, the algebra
$A_\sigma$ is $(H,\,r)$-Azumaya if and only if $\theta_\sigma$ is invertible. 
\end{theorem} 
\pf Let $\theta_\sigma$ be invertible. We shall see that, for every
$\eta\in A_{\sigma}^*$ and every $m\in A_{\sigma}$ the endomorphism of
$A_{\sigma}$ given by $h\mapsto \langle \eta,\,h\rangle m$ belongs to
the image of $F$. Let $h\in H$ be such that $\eta=w(h)$ and let us consider the following element
of $A_{\sigma}\#\overline{A_{\sigma}}$: 
$$\Gamma=\sum m\cdot S_2(S(\h2))\cdot
S_1(\theta_{\sigma}^{-1}((w(\h1))_{(2)}))
\#\theta_{\sigma}^{-1}((w(\h1))_{(1)}).$$ Then for every $l\in A_{\sigma}$
$$
\begin{array}{rl}
F(\Gamma)(l)&=\sum\langle \theta_{\sigma}
(\theta_{\sigma}^{-1}((w(\h1)_{(1)}),\,
\l1\rangle \,m\cdot S_2(S(\h2))\cdot \\
&\phantom{=\sum} S_1(\theta_{\sigma}^{-1}((w(\h1))_{(3)}))\cdot\theta_{\sigma}^{-1}((w(\h1))_{(2)})\cdot\l2\\
&=\sum \langle (w(\h1))_{(1)},\,\l1\rangle \, m\cdot S_2(S(\h2))\cdot 
\varepsilon((w(\h1))_{(2)})\l2\\
&=\sum\langle w(\h1),\,\l1\rangle \,m\cdot S_1(S(\h2))\cdot \l2
\end{array}
$$
Applying (\ref{vu}) we have
$$
\begin{array}{rl}
F(\Gamma)(l)&=\sum\langle w(\h1),\,l\rangle \,m\cdot S_2(S(\h3))\cdot S(\h2)\\
&=\sum\langle w(\h1),\,l\rangle \,m\cdot S_2((S(\h2))_{(1)})\cdot 
(S(\h2))_{(2)}\\
&=\sum \langle w(\h1),\,l\rangle  \varepsilon(S(\h2))\,m=
\langle w(h),\,l\rangle m
\end{array}
$$
where we used (\ref{vare}). Hence $F$ is surjective.
%
Similarly, let $(\theta_\sigma^{-1})^*$ be 
the dual map of $\theta_\sigma^{-1}$ with respect to the
non-degenerate pairing $\langle\phantom{-},\phantom{-}\rangle$. 
The map $(\theta_\sigma^{-1})^*$ is a well-defined Hopf algebra map
$(_\sigma H_{\sigma^{-1}})^{*, cop}\longrightarrow\,_\sigma
H_{\sigma^{-1}}$ and it is bijective if $\theta_\sigma$ is so. 
Let $h$ and $m$ be as before and let $\Gamma'$ be the following
element of  $A_{\sigma}\#\overline{A_{\sigma}}$: 
$$\Gamma'=\sum (\theta^{-1}_{\sigma})^*((w(\h1))_{(2)})\#
 S_2((\theta^{-1}_{\sigma})^*((w(\h1))_{(1)}))\cdot S_1(S(\h2))\cdot
m.$$ Then we have: 
$$
\begin{array}{rl}
G(\Gamma')(l)&=
\sum \langle \theta_{\sigma}(\l1),\,((\theta^{-1}_{\sigma})^*
((w(\h1))_{(2)}))_{(1)}\rangle \,\l2\cdot\\
&\phantom{,}((\theta^{-1}_{\sigma})^*((w(\h1))_{(2)}))_{(2)}\cdot S_2
((\theta^{-1}_{\sigma})^*((w(\h1))_{(1)}))\cdot S_1(S(\h2))\cdot m \\
&=\sum \langle \theta_{\sigma}(\l1),\,(\theta^{-1}_{\sigma})^*
(((w(\h1))_{(2)})_{(2)}))\rangle \,\l2\cdot\\
&\phantom{,}(\theta^{-1}_{\sigma})^*(((w(\h1))_{(2)})_{(1)})\cdot S_2
((\theta^{-1}_{\sigma})^*((w(\h1))_{(1)})\cdot S_1(S(\h2))\cdot m \\
&=\sum \langle \theta_{\sigma}(\l1),\,(\theta^{-1}_{\sigma})^*
((w(\h1))_{(3)}))\rangle \,\l2\cdot\\
&\phantom{,}(\theta^{-1}_{\sigma})^*((w(\h1))_{(2)})\cdot S_2
((\theta^{-1}_{\sigma})^*((w(\h1))_{(1)})\cdot S_1(S(\h2))\cdot m\\
& =\sum \langle \theta_{\sigma}(\l1),\,(\theta^{-1}_{\sigma})^*
((w(\h1))_{(2)})\rangle \,\l2\cdot((\theta^{-1}_{\sigma})^*((w(\h1))_{(1)})_{(2)}\cdot\\
&\phantom{,}\cdot S_2
((\theta^{-1}_{\sigma})^*(((w(\h1))_{(1)})_{(1)})\cdot S_1
(S(\h2))\cdot m \\
&=\sum \langle \theta_{\sigma}(\l1),\,(\theta^{-1}_{\sigma})^*
((w(\h1))_{(2)}))\rangle \,\l2\cdot((\theta^{-1}_{\sigma})^*((w(\h1))_{(1)}))_{(1)}\cdot\\
&\phantom{,}S_2
((\theta^{-1}_{\sigma})^*(((w(\h1))_{(1)}))_{(2)})\cdot S_1
(S(\h2))\cdot m\\
&=\sum \langle w(\h1),\,\l1\rangle \,\l2\cdot S_1(S(\h2))\cdot m
\end{array}
$$
where for the last equality we used (\ref{vare3}) applied to
$(\theta_\sigma^{-1})^*((w(\h1))_{(1)})$. By (\ref{vu}) we get:
$$
\begin{array}{ll}
G(\Gamma')(l)&=\sum \langle w(\h1),\,l\rangle \,S(\h2)\cdot 
S_1(S(\h3))\cdot m\\
&=\sum \langle w(\h1),\,l\rangle \,(S(\h2))_{(2)}\cdot S_1((S(\h2))_{(1)})\cdot m\\
&=\sum \langle w(h),\,l\rangle \,m
\end{array}
$$
where the last equality follows from (\ref{vare}).
Therefore, if $\theta_{\sigma}$ is bijective then $F$ and $G$ are
surjective, hence bijective. 

\smallskip

Let us now assume that $F$ is bijective. We will show that
$\theta_\sigma$ is surjective. We recall that if $A$ is a right
$H$-comodule, $End(A)$ with comodule structure
given by  (\ref{endo}) 
is isomorphic, as right comodule, to 
$A\otimes A^*$ with comodule structure on $A^*$ given by 
$$\rho(\xi)(a)=\sum \langle\xi,\,\a1\rangle\,S^{-1}(\a2)$$
for every $\xi\in A^*$ and every $a\in A$. An isomorphism $\Phi$ 
is given by
$\Phi(a\otimes \xi)(b)=\langle \xi,\,b\rangle\,a$ for every
$a,\,b\in A$ and every $\xi\in A^*$. 

Let us consider $\eta\in H^*$. We shall show that $\eta$ belongs to the
image of $\theta_\sigma$. We know that $\eta=w(h)$ for some
$h$. For $A=A_\sigma\cong H$ as comodule,
$\Phi(1\otimes\eta)$ belongs to the image of $F$, hence, there exists a 
$\Gamma=\sum a_i\# b_i\in A_\sigma\#\overline{A_\sigma}$ for which
$F(\Gamma)=1\otimes\eta=1\otimes w(h)$.
Since $F$ is a comodule map, there holds:
$$F(\sum a_{i1}\# b_{i1})\otimes b_{i2}a_{i2}=\sum 1\otimes\eta_{(1)}\otimes\eta_{(2)}$$
that is, for every $k\in H$
$$
\begin{array}{l}
\sum\langle\theta_\sigma(b_{i(1)}),\,\k1\rangle\,a_{i(1)}\cdot b_{i(2)}\cdot \k2\otimes b_{i(3)}a_{i(2)}=\sum\langle\eta,\,\k1\rangle\,S^{-1}(\k2)=\\
\sum\langle w(h),\,\k1\rangle\,S^{-1}(\k2)=\sum \langle w(\h1),\,k\rangle\,\h2
\end{array}
$$
where the last equality follows from (\ref{pallino}).
Applying the linear operator $(S\otimes\id)\tau$ to the
last and the first term of the above chain of equalities we obtain: 
$$
\begin{array}{rl}
\sum\!\!\! &\langle w(\h1),k\rangle S(\h2)\otimes 1\\
&=\sum S(a_{i(2)})S(b_{i(3)})\otimes \langle\theta_\sigma(b_{i(1)}),\,
\k1\rangle\,a_{i(1)}\cdot b_{i(2)}\cdot \k2\\
&=\sum \sigma^{-1}((S(a_{i(2)}))_{(1)}\otimes(S(b_{i(3)}))_{(1)})
(S(b_{i(3)}))_{(2)}\cdot (S(a_{i(2)}))_{(2)}\otimes\\
&\phantom{=}\langle\theta_\sigma(b_{i(1)}),\,\k1\rangle\,a_{i(1)}\cdot b_{i(2)}\cdot \k2\\
&=\sum \sigma^{-1}(S(a_{i(3)})\otimes S(b_{i(4)}))S(b_{i(3)})\cdot 
S(a_{i(2)})\otimes\\ 
&\phantom{=}\langle\theta_\sigma(b_{i(1)}),\,\k1\rangle\,a_{i(1)}\cdot
b_{i(2)}\cdot \k2.
\end{array}
$$
Applying the product in $A_{\sigma}$ on the first and the last term of
the above chain of equalities we obtain: 
$$
\begin{array}{rl}
\sum \langle w(\h1),k\rangle S(\h2)&=\sum \sigma^{-1}(S(a_{i(3)})\otimes 
S(b_{i(4)}))S(b_{i(3)})\cdot S(a_{i(2)})\cdot \\
&\phantom{=}\langle\theta_\sigma(b_{i(1)}),\,\k1\rangle\,a_{i(1)}\cdot b_{i(2)}\cdot \k2.
\end{array}
$$
A direct computation yields, for every $l\in H$:
$$\sum S(\l2)\cdot\l1=\sum\sigma(\l1\otimes S(\l4))\l2\,S(\l3)=
\sum\sigma(\l1\otimes S(\l2)).$$
Using this formula in the previous equality we have:
$$ 
\begin{array}{ll}
\sum \langle w(\h1),\,k\rangle S(\h2)&=\sum \sigma^{-1}(S(a_{i(3)})\otimes S(b_{i(4)}))\sigma(b_{i(2)}\otimes\\ 
&\phantom{=}S(b_{i(3)}))\sigma (a_{i(1)}\otimes S(a_{i(2)}))
\langle\theta_\sigma(b_{i(1)}),\,\k1\rangle\,\k2. 
\end{array}
$$
Applying $\varepsilon$ on both sides and observing that the equality holds for every
$k$ yields:
$$
\eta=\theta_\sigma\left(\sum\sigma^{-1}(S(a_{i(3)})
\otimes
S(b_{i(4)}))\sigma(b_{i(2)}\otimes 
S(b_{i(3)}))\sigma (a_{i(1)}\otimes
S(a_{i(2)}))(b_{i(1)})\right).
$$
Hence, $\eta\in
Im(\theta_\sigma)$ for every $\eta\in H^*$. 

Let us now suppose that $G$ is surjective. In a similar fashion we
shall prove that $\theta_\sigma^*$ is surjective. 
The right $H$-comodule $End(A)^{op}$
with comodule structure given by (\ref{endop})
is isomorphic to the right $H$-comodule  
$A^*\otimes A$ with comodule structure on $A^*$ given by
$$\rho(\xi)(a)=\sum \langle\xi,\,\a1\rangle\,S(\a2)$$
for every $\xi\in A^*$ and every $a\in A$. An isomorphism is given by
$\Psi(\xi\otimes a)(b)=\langle \xi,\,b\rangle\,a$ for every
$a,\,b\in A$ and every $\xi\in A^*$. 

Let $u(m)=S^*(v(m))$ for every $m\in H$, let $\eta\in H^*$ and let $h$
be such that 
$\eta=u(h)$. As before  $\Psi(\eta\otimes 1)=G(\Gamma')$ for some 
$\Gamma'=\sum c_i\# d_i\in A_\sigma\#\overline{A_\sigma}$.
Since $G$ is a comodule map, there holds:
$$G(\sum c_{i(1)}\# d_{i(1)})\otimes d_{i(2)}c_{i(2)}=\sum
\eta_{(1)}\otimes1\otimes\eta_{(2)}$$ 
that is, for every $k\in H$:
$$
\begin{array}{c}
\sum\langle\theta_\sigma(\k1),\,c_{i(1)}\rangle\,\k2\cdot
c_{i(2)}\cdot d_{i(1)}\otimes d_{i(2)}c_{i(3)}=
\sum\langle\eta,\,\k1\rangle\,1\otimes S(\k2)\\
=\sum\langle u(h),\,\k1\rangle\,1\otimes\ S(\k2)=\sum \langle v(h),
(S(k))_{(2)}\rangle 1\otimes (S(k))_{(1)}\\
=\sum \langle u(\h1),\,k\rangle\,1\otimes \h2.
\end{array}
$$
Applying the linear operator $(\id\otimes S^{-1})$ to the first and
the last term of the above chain of equalities we obtain: 
$$
\begin{array}{ll}
\sum \langle u(\h1),\,k\rangle\,1&\otimes S^{-1}(\h2)\\
&=\sum \langle\theta_\sigma(\k1),\,c_{i(1)}\rangle\,\k2\cdot
c_{i(2)}\cdot d_{i(1)}\otimes S^{-1}(c_{i(3)})S^{-1}(d_{i(2)})\\
&=\sum \langle\theta_\sigma(\k1),\,c_{i(1)}\rangle\,\k2\cdot
c_{i(2)}\cdot d_{i(1)}\otimes\\
&\phantom{=} S^{-1}(d_{i(2)})\cdot
S^{-1}(c_{i(3)})\sigma^{-1}(S^{-1}(c_{i(4)})\otimes S^{-1}(d_{i(3)})). 
\end{array}
$$
Applying the product in $A_{\sigma}$ on the first and the last term of
the above chain of equalities and using the formula: 
$$
\sum \l1\cdot
S^{-1}(\l2)=\sum\sigma(S^{-1}(\l4)\otimes\l1)S^{-1}(\l3)\l2=
\sum\sigma(S^{-1}(\l2)\otimes\l1)  
$$ 
for every $l\in H$ we obtain:
$$
\begin{array}{ll}
\sum \langle u(\h1),k\rangle&S^{-1}(\h2)=
\sum \langle\theta_\sigma(\k1),\,c_{i(1)}\rangle\, \sigma(S^{-1}(d_{i(2)})\otimes d_{i(1)})\\
&\sigma(S^{-1}(c_{i(3)})\otimes
c_{i(2)})\sigma^{-1}(S^{-1}(c_{i(4)})\otimes  
S^{-1}(d_{i(3)}))\k2.
\end{array}
$$
Applying $\varepsilon$ on both sides and observing that equality holds
for every $k$ yields $\eta= \theta_{\sigma}^*(z)$ with
$$
z=\sum \sigma(S^{-1}(d_{i(2)})\otimes
d_{i(1)})\sigma(S^{-1}(c_{i(3)})\otimes
c_{i(2)})\sigma^{-1}(S^{-1}(c_{i(4)})\otimes  
S^{-1}(d_{i(3)}))c_{i(1)}
$$
whence the proof.\hfill$\Box$

\begin{corollary}Let $H$ be a finite-dimensional  dual quasitriangular
Hopf algebra with universal $r$-form $r$. Then 
the Hopf algebra $H^{op}$ is $(H,\,r)$-Azumaya if and only if
$\theta_r$ is bijective.\hfill$\Box$ 
\end{corollary}

\begin{remark}{\rm Let us observe that, as a consequence of the proof of
 Theorem \ref{analogo}, $F$ is bijective if and only if $G$ is so. }
\end{remark}

\begin{remark}{\rm The condition in Theorem \ref{analogo} that the Hopf algebra is over a field
  could be relaxed. Indeed the proof would work by localization, just as in
  \cite[Theorem 12.4.5]{caen}, for any faithfully projective Hopf
  algebra with bijective antipode over a commutative ring.}   
\end{remark}


\begin{remark}{\rm Let $(H,\,r)$ be dual quasitriangular. Then
$(H^{op},\,r\tau)$ is again dual quasitriangular. Theorem
\ref{analogo} for $H^{op}$ states that if $\sigma$ is a $2$-cocycle
for $H^{op}$ 
then $_{\sigma\tau}H$ is $(H^{op},\,r\tau)$-Azumaya if and only if the
map: 
$$\theta_\sigma\colon H\longrightarrow (H^{op})^*$$
$$h\mapsto r_\sigma(h\otimes\,-)$$ 
is an isomorphism. If $H$ is commutative and cocommutative, we recover 
\cite[Theorem 12.4.5]{caen} with $f=\sigma\tau$ and 
$\theta\colon H\longrightarrow H^*$ given by
$\theta(h)=r\tau(h\otimes\,-).$ Indeed, the action of $H$ determined
by the pairing $r\tau$ is 
$$h\lu k=\sum\k1\,r(\k2\otimes h)=\sum \k1\,\langle\theta(\k2),\,h\rangle.$$
Let us observe that our $\theta_\sigma$ is
Caenepeel's $(\theta*d)^*$.}\end{remark}

%

\subsection{The dual picture}

We would like to outline briefly the dual picture, i.e., the analysis
of $(H,\,R)$-Azumaya algebras for a quasitriangular Hopf algebra $H$
with $R$-matrix $R$. 
In order to fix notation we recall well-known facts about isomorphisms
of Brauer groups and the standard equivalence $\D$ between the 
category $\;_H\!\M$ and the category $\M^{H^*}$ (see for instance
\cite[Lemma 1.6.4]{M}). Let $\otimes^{op}$ denote the functor obtained
from $\otimes$ by reversing the order of the tensorands. 
The functor $\D$, together with the natural transformations
$\tau_{UV}\colon \D(U)\otimes^{op}\D(V)\to \D(U\otimes V)$ for every
pair of objects $U,V$ in $_H\M$, and with $\id_U\colon \f\to\D(\f)=\f$, define an
equivalence of monoidal categories between $(_H\!\M,\otimes, \f)$
and $(\M^{H^*},\otimes^{op}, \f)$. If $H$ is quasitriangular with
$R$-matrix $R=\sum R^{(1)}\otimes R^{(2)}$ then $H^*$ is dual
quasitriangular with universal 
$r$-form $R$, viewed as an element of $(H\otimes H)^{**}$. 
The functor $\D$ together with $\tau$ and $\id$
define an equivalence of braided  monoidal categories between
$(\M^{H^*},\otimes,\f)$ and $(_H\M,\otimes^{op},\f)$.
Here the braiding in $(\M^{H^*},\otimes,\f)$ is given by
$(u\otimes v)\mapsto\sum R^{(2)}.v\otimes R^{(1)}.u$, the
braiding in $(\M^{H^*},\otimes, \f)$ is given
by  $m\otimes n\mapsto\sum\n0\otimes\m0 \langle \n1,R^{(1)}\rangle\,\langle
\m1, R^{(2)}\rangle$ and the braiding in  $(\M^{H^*},\otimes^{op}, \f)$ is induced by
the braiding in  $(\M^{H^*},\otimes, \f)$. The reversed equivalence
induces an isomorphism 
$$BM(\f, H, R)\mapsto BC(\f, H^*, R)^{op}\cong
BC(\f, H^*, R)$$ where the class of $A$ with given $H$-module structure 
is mapped to the class of $A^{op}$ with right $H^{*,op}$-comodule
structure determined by the functor $\D$.

The dual version of Theorem \ref{analogo} reads:
\begin{corollary}\label{duanalogo}
Let $H$ be a finite-dimensional quasitriangular  Hopf algebra with
$R$-matrix $R$. Let  $C=\sum C^{(1)}\otimes C^{(2)}\in H\otimes H$ 
be a cocycle 
for $H^*$ and let $R_C=(\tau C)R C^{-1}$. Then $\f\#_C H^*=\,_C H^*$
with $H$-action: $h\lu 
f=\sum f_{(1)}\langle f_{(2)}, h\rangle$ is $(H, R)$-Azumaya if and
only if the map: $\theta\colon H^{*,op}\longrightarrow H$ given by 
$\theta(f)=\sum R_C^{(1)}\langle f,\, R_C^{(2)}\rangle$ 
is an isomorphism.\hfill$\Box$  
\end{corollary}  
In particular
\begin{corollary}Let $H$ be a finite-dimensional quasitriangular  Hopf
  algebra with 
$R$-matrix $R$. Then $H^*$ with $H$-action: $h\lu
f=\sum f_{(1)}\langle f_{(2)}, h\rangle$ is $(H, R)$-Azumaya if and
only if the map: $\theta\colon H^{*,op}\longrightarrow H$ given by 
$\theta(f)=\sum R^{(1)}\langle f,\, R^{(2)}\rangle$ 
is an isomorphism.\hfill$\Box$ 
\end{corollary}
  
\section{An Example: $E(n)$}

Let $char(\f)\neq2$, let $n\ge1$ and let $E(n)$ denote the Hopf
  algebra generated 
  by $c$ and $x_i$ for $1\le i\le n$ with relations: 
$$c^2=1,\quad cx_i+x_ic=0;\quad x_ix_j+x_jx_i=0,\quad x_i^2=0$$ coproduct: 
$$\Delta(c)=c\otimes c;\quad \Delta(x_i)=1\otimes x_i+x_i\otimes c$$
and antipode $S(c)=c$ and $S(x_j)=cx_j$. 
The Hopf algebra $E(n)$ is quasitriangular, isomorphic to its
opposite and self-dual. Its $R$-matrices were classified in \cite{PVO1} and they are 
parametrized by matrices in $M_n(\f)$. By self-duality, the universal $r$-forms 
are parametrized by elements of $M_n(\f)$ and they are given as
follows: for a matrix $A=(a_{ij}) \in M_n(\f)$ and for $s$-tuples
$P,F$ of increasing elements in $\{1,\ldots,n\}$ we define
$|P|=|F|=s$ and $x_P$ as the product of the $x_j$'s whose index
belongs to $P$, taken in increasing order. Any bijective map
$\eta:F \rightarrow F$ may be identified with an element of the
symmetric group $S_s$. Let $sign(\eta)$ denote the signature of
$\eta$. If $P=\emptyset$ then we take $F=\emptyset$ and
$sign(\eta)=1$. Finally, by $a_{P,\eta(F)}$ we denote the product
$a_{p_1,f_{\eta(1)}}\cdots\,a_{p_s,f_{\eta(s)}}$. For
$P=\emptyset$ we define $a_{P,\eta(F)}:=1$. Then the universal
$r$-form corresponding to the matrix $A$ is:
$$\begin{array}{l}
r_A= \sum_P (-1)^{{|P|(|P|-1)}\over2}\sum_{F, |F|=|P|, \eta\in
S_{|P|}}sign(\eta)a_{P,\eta(F)}\bigl((x_P)^*\otimes (x_F)^*  \\
  +(cx_{P})^*\otimes (x_F)^*+(-1)^{|P|}(x_P)^*\otimes (cx_F)^*-(-1)^{|P|}
 (cx_{P})^*\otimes (cx_F)^*\bigr).
\end{array}$$
In particular, $r_A(x_i\otimes x_j)=a_{ij}.$

\smallskip

The $E(n)$-cleft extensions of $\f$ up to equivalence were classified 
in \cite{PVO2}. They are  parametrized by an invertible scalar
$\alpha$, a vector $\gamma\in\f^n$ and a lower triangular $n\times n$
matrix $\Lambda=(\lambda_{ij})$. 
On the generators the corresponding cocycle $\sigma=\sigma(\alpha,
\gamma, \Lambda)$ has values:
$$\sigma(c\otimes c)=\alpha;\quad \sigma(x_i\otimes
c)=\gamma_i;\quad\sigma(c\otimes x_i)=0;\qquad\sigma(x_i\otimes x_j)=\lambda_{ij}.$$
The cleft extension corresponding to
the cocycle $\sigma$ is the
generalized Clifford algebra $Cl(\alpha,\gamma,\lambda)$ 
with generators $u$ and $v_i$, for $i=1,\,\ldots n$, relations
\begin{eqnarray}\label{privileged}
& u^2=1,\ uv_i+v_iu=\gamma_i, \ & v_j^2=\lambda_{jj}, \
v_iv_j+v_jv_i=\lambda_{ij} \quad \mbox{ for $i\neq j$}\quad
\end{eqnarray}
and with comodule algebra structure given by:
\begin{equation}\label{privileged2}
\rho(u)=u\otimes c,\quad \rho(v_j)=1\otimes x_j+v_j\otimes c.
\end{equation}
It is clear that $_{\sigma\tau}
E(n)^{op}=Cl(\alpha,\gamma,\Lambda)^{op}\cong
Cl(\alpha,\gamma,\Lambda)$. We shall apply Theorem \ref{analogo} to
$Cl(\alpha,\gamma,\lambda)$ and reduce the question on when this
comodule algebra is $(E(n),r_A)$-Azumaya to a simple linear algebra
problem. Since 
$E(n)=(_\sigma E(n)_{\sigma^{-1}})^{op}$ as coalgebras, by
\cite[Proposition 2.4.2]{hr} the coalgebra map
$\theta_\sigma\colon (_\sigma E(n)_{\sigma^{-1}})^{op}\to\;
(_\sigma E(n)_{\sigma^{-1}})^{*}$ is injective if and only if its
restriction to the 
span $W$ of $1,c$  the $x_j's$ and the $cx_j$'s is injective. 
For every $h\in E(n)$ let us denote the corresponding element in the twisted
Hopf algebra $_\sigma E(n)_{\sigma^{-1}}$ by $\overline{h}$. 
If for every $\overline{h}\in W$ the restriction of the functional
$\theta_\sigma(\overline{h})$ to $W$ is not identically zero, then
$\theta_\sigma$ is injective. Let us now assume that
there exists an element $\overline{h}\in W$ such that 
$\theta_\sigma(\overline{h})(\overline{k})=0$ for every
$\overline{k}\in W$. Let 
$$\overline{h}=e+f\overline{c}+\sum s_i\overline{x}_i+\sum
t_i\overline{cx_i}.$$

By the description of the cocycles in \cite{PVO2} we have:
$$\sigma(cx_i\otimes c)=\gamma_i;\quad \sigma(cx_i\otimes x_j)=\lambda_{ij};$$
 $$\sigma(c\otimes cx_i)=0;\quad \sigma(x_i\otimes
cx_j)=-\lambda_{ij};\quad\sigma(cx_i\otimes
cx_j)=-\alpha\lambda_{ij}$$
and
$$\begin{array}{lll}\\
\sigma^{-1}(c\otimes c)=\alpha^{-1};&\sigma^{-1}(c\otimes
x_j)=0;&\sigma^{-1}(x_j\otimes c)=-\alpha^{-1}\gamma_j;\\
\sigma^{-1}(c\otimes cx_j)=0;&\sigma^{-1}(cx_j\otimes
c)=-\alpha^{-1}\gamma_j;&\sigma^{-1}(x_i\otimes x_j)=-\alpha^{-1}\lambda_{ij};\\
\sigma^{-1}(cx_i\otimes x_j)=-\lambda_{ij};&\sigma^{-1}(x_i\otimes
cx_j)=\lambda_{ij};&\sigma^{-1}(cx_i\otimes cx_j)=\lambda_{ij}.
\end{array}
$$
A direct computation shows that, for the universal $r$-form $r_A$ and
the above cocycle $\sigma$ one has: 
$$\begin{array}{lll}\\
r_{A,\sigma}(x_j\otimes c)=-\alpha^{-1}\gamma_j,&r_{A,\sigma}(c\otimes
x_j)=-\alpha^{-1}\gamma_j,&r_{A,\sigma}(c\otimes cx_j)=\gamma_j,\\
r_{A,\sigma}(cx_i\otimes cx_j)=\alpha b_{ij,}&r_{A,\sigma}(x_i\otimes x_j)=\alpha^{-1}b_{ij},&r_{A,\sigma}(cx_j\otimes
c)=\gamma_j,\\
r_{A,\sigma}(cx_i\otimes x_j)=b_{ij}+\alpha^{-1}\gamma_i\gamma_j,&r_{A,\sigma}(x_i\otimes
cx_j)=-b_{ij},&r_{A,\sigma}(c\otimes c)=-1
\end{array}
$$
where $B$ is the $n\times n$ matrix whose $(i,j)$-entry is
$b_{ij}=a_{ij}-\lambda_{ij}-\lambda_{ji}$. 

\smallskip

We have:
$$
\begin{array}{ll}
\theta_\sigma(\overline{h})(1)&=e+f=0\\
\theta_{A,\sigma}(\overline{h})(\overline{c})&=e+fr_{A,\sigma}(c\otimes
 c)+\sum_js_jr_{A,\sigma}(c\otimes x_j)+\sum_jt_jr_{A,\sigma}
(c\otimes cx_j)\\
&=e-f+\sum(t_j-\alpha^{-1}\sum s_j)\gamma_j=0\\
\theta_{A,\sigma}(\overline{h})(\overline{x_i})&=-\alpha^{-1}f\gamma_i+\sum b_{ij}(\alpha^{-1}s_j-t_j)=0\\
\theta_{A,\sigma}(\overline{h})(\overline{cx_i})&=f\gamma_i+\sum
b_{ij}(s_j+\alpha t_j)+
\alpha^{-1}\gamma_i\sum \gamma_js_j=0.
\end{array}
$$
This is possibe if and only if there exists a non-trivial solution
$(y,\underline{x},\underline{z})$ with $y\in\f$ and $\underline{x},
\underline{z}\in\f^n$ to the system:
$$
\left\{
\begin{array}{ll}
(\underline{z}-\alpha^{-1}\underline{x})\bullet\gamma&=2y;\\
B(\alpha\underline{z}-\underline{x})&=-y\gamma;\\
B(\alpha\underline{z}+\underline{x})&=-(y+\alpha^{-1}\underline{x}\bullet\gamma)\gamma
\end{array}\right.
$$
where $\bullet$ denotes the usual dot product in $\f^n$. This system
is equivalent to:
$$
\left\{
\begin{array}{ll}
\underline{z}\bullet\gamma&=2y+\alpha^{-1}\underline{x}\bullet\gamma;\\
2\alpha B\underline{z}&=-(2y+\alpha^{-1}\underline{x}\bullet\gamma)\gamma;\\
2B\underline{x}&=-(\alpha^{-1}\underline{x}\bullet\gamma)\gamma
\end{array}\right.
$$
which is equivalent to:
$$
\left\{
\begin{array}{ll}
\underline{z}\bullet\gamma&=2y+\alpha^{-1}\underline{x}\bullet\gamma;\\
2\alpha B\underline{z}&=-(\underline{z}\bullet\gamma)\gamma;\\
2\alpha B\underline{x}&=-(\underline{x}\bullet\gamma)\gamma.
\end{array}\right.
$$
If the third equation admits a non-trivial solution $\underline{s}$,
 we may take $\underline{z}=0$ and 
$2f=2y=-\alpha^{-1}\underline{s}\bullet\gamma$ 
and the system admits a
non-trivial solution. If the third equation does not admit a
non-trivial solution then the same holds for the second equation
forcing $\underline{x}=\underline{0}=\underline{z}$ and $y=0$. 
Hence, if such an $\overline{h}$ exists, we may assume that it is $(1,c)$-skew-primitive. 

Let us observe that, due to the particular coalgebra structure of
$E(n)$, the elements $\overline{c}$ and $\overline{x}_j$ for
$j=1,\ldots,n$ are algebra generators in $_\sigma E(n)_{\sigma^{-1}}$. Indeed,
the elements $\overline{c^ax}_P$ with $a=0,1$ and
$P\subset\{1,\ldots,n\}$  span 
$_\sigma E(n)_{\sigma^{-1}}$. 
One can prove by induction on $|P|=m$ that 
$\overline{cx}_P$ lies in the span of $\overline{c}\cdot_\sigma\overline{x}_P$ 
and of terms of the form 
$\overline{c}^b\cdot_\sigma\overline{x_{j_1}}\cdot_\sigma\cdots\cdot_\sigma\overline{x_{j_k}}$ with
$J$ strictly contained in $P$ and $b=0,1$. Similarly, it can be proved
bu induction on $|P|=m$ that  
$\overline{x}_P$ lies in the span of
$\overline{x}_{p_1}\cdot_\sigma\cdots\cdot_\sigma\overline{x}_{p_m}$ and 
$\overline{c}^b\cdot_\sigma\overline{x_{j_1}}\cdot_\sigma\cdots\cdot_\sigma\overline{x_{j_k}}$ with
$J$ strictly contained in $P$ and $b=0,1$. Then 
every element $\overline l$  of $_\sigma E(n)_{\sigma^{-1}}$ is
spanned by a
product 
$\overline{l}=\overline{k}_1\cdot_\sigma\cdots\cdot_\sigma\overline{k}_m$ with $k_j\in W$. 
Thus 
$$
\begin{array}{rl}
(\theta_\sigma(\overline{h}))(\overline{l})&=\theta_\sigma(\overline{h})(\overline{k}_1\cdot_\sigma\cdots\cdot_\sigma\overline{k}_m)\\
&=\sum \prod \theta_\sigma(\overline{h}_{(i)})(\overline{k}_i)=0
\end{array}
$$
because in each summand
we have $h_{(i)}=h$ for some $i$, since we have assumed $h$ to be
skew-primitive. 
Hence if $\theta_\sigma(\overline{h})$ is identically zero on $W$, it
is zero and $\theta_\sigma$ is not injective. 
If we denote by $\Gamma$ is the $n\times n$ matrix whose $(i,j)$-entry
is $\Gamma_{ij}=\gamma_i\gamma_j$, the previous discussion shows that
$\theta_\sigma$ is injective if and only
if$2\alpha B\underline{x}=-(\underline{x}\bullet\gamma)\gamma$, and
therefore  
$$(2\alpha B+\Gamma)\underline{z}=\underline{0}$$
admits no non-trivial solutions. We have reached the following result:
\begin{proposition}\label{wk}The cleft extension $Cl(\alpha,\gamma, \Lambda)^{op}$ is
  $(E(n), r_A)$-Azumaya if and only if
  $\det(2\alpha(A-\Lambda-\Lambda^t)+\Gamma)\neq0$. \hfill$\Box$
\end{proposition}
It was proved in \cite[Lemma 2.3]{GioJuan3} that when $\gamma=0$ and
$\alpha=1$ the corresponding cocycle $\sigma$ is lazy, i.e., the
product in the twisted Hopf algebra $_\sigma E(n)_\sigma^{-1}$
coincides with the product in $E(n)$. In this case the map
$\theta_\sigma$ is again a Hopf algebra map $E(n)^{op}\to E(n)^*$ and
Proposition \ref{wk} states that $Cl(1,0, \Lambda)^{op}$ is
  $(E(n), r_A)$-Azumaya if and only if
  $\det(A-\Lambda-\Lambda^t)\neq0$. In particular, when $\sigma$ is
trivial we have:
\begin{corollary}The comodule algebra $E(n)^{op}$ is $(E(n),
  r_A)$-Azumaya if and only if $\det(A)\neq0$. 
\end{corollary}
In the computation of the Brauer groups of $H_4$, $E(n)$ and $H_\nu$
a special r{o}le is played by those universal $r$-forms that are non-zero
only on the group algebra of the grouplike elements. For $E(n)$ this
is $r_0$. In this case we have:
\begin{corollary}\label{cotri}The comodule algebra $Cl(1,0,\Lambda)$ is $(E(n),
  r_0)$-Azumaya if and only if $\det(\Lambda+\Lambda^t)\neq0$.
\end{corollary}

It was shown in \cite{GioJuan3} that $BM(\f, E(n), R_0)$ is
isomorphic to the direct product of the Brauer-Wall group
$BW(\f)$ of the field
$\f$ and the group $Sym_n(\f)$ 
of $n\times n$ symmetric matrices with coefficients
in $\f$, (represented by special cocycles cohomologous to those in
\cite{PVO2}). On the other hand, the map $\phi\colon
E(n)\to E(n)^*$ with $\phi(c)=1^*-c^*$ and
$\phi(x_j)=x_j^*+(cx_j)^*$ for $j=1,\ldots,n$ defines a Hopf algebra
isomorphism. Therefore the  pull-back along $\phi$ yields an
isomorphism 
$$BM(\f, E(n)^*,r_0)=BM(\f, E(n)^*, (\phi\otimes\phi)(R_0))\cong
BM(\f, E(n), R_0).$$  
Since $BC(\f, E(n), r_0)\cong
BM(\f, E(n)^*, r_0)$, we may identify $BC(\f, E(n), r_0)$ with $BM(\f,
E(n), R_0)$. The class of $Cl(\alpha,\gamma,
\Lambda)^{op}$ in $BC(\f, E(n), r_0)$ 
described above corresponds to the class of the algebra  
$Cl(\alpha,\gamma,
\Lambda)$ with action: 
$$
\begin{array}{ll}
c\lu u=u\langle\phi(c), c\rangle=-u;& 
c\lu v_i=v_i\langle\phi(c), c\rangle+1\langle\phi(c),
x_i\rangle=-v_i;\\
x_i\lu u=u\langle\phi(x_i), c\rangle=0;&x_j\lu
v_i=v_i\langle\phi(x_j),c\rangle+1\langle\phi(x_j), 
x_i\rangle=\delta_{ij}.
\end{array}
$$
We end this section describing the decomposition of the class
represented by $C(\alpha,\gamma,\Lambda)$ as
a product of an element in $BW(\f)$ and an element in $Sym_n(\f)$. 
We recall that the product $\#$ corresponding to $R_0$ and $r_0$ is
just the ${\mathbb Z}_2$-graded tensor product, where the grading is
induced by the eigenspace decomposition with respect to the action of
$c$.

\smallskip 

Let us observe that, taking $z_i=v_i-\frac{\gamma_i}{2\alpha}u$, the
algebra $Cl(\alpha,\gamma,
\Lambda)$ is isomorphic, as an $E(n)$-module algebra, to the algebra
with generators $u$ and $z_1,\ldots,z_n$, relations:
$$
\begin{array}{cc}  
uz_i+z_iu=0;&z_iz_j+z_jz_i=2\left(\frac{1}{4\alpha}\right)\left(2\alpha(\Lambda_{ij}+\Lambda_{ji})-\Gamma_{ij}\right);\\
\\
u^2=\alpha;&z_i^2=\left(\frac{1}{4\alpha}\right)\left(4\alpha\Lambda_{ii}-\Gamma_{ij}\right) 
\end{array}
 $$ and with action: 
$$
\begin{array}{cc}
c\lu u=-u;& 
c\lu z_i=-z_i;\\
x_i\lu u=0;&x_j\lu
z_i=\delta_{ij}.
\end{array}
$$ 
Hence, we may always reduce
to the case that the cleft extension is associated to a cocycle $\omega$ with
$\omega(x_i\otimes x_j)=l_{ij}$ with $L$ a symmetric matrix,
and $\omega(c\otimes x_j)=\omega(x_j\otimes c)=0$ (see also \cite[\S
  2]{GioJuan3}). We shall denote such a module algebra by $C(\alpha,
L)$. We observe that the map $u\mapsto tu$ gives an isomorphism $C(\alpha,
L)\cong C(t^2\alpha, L)$.
The algebra $C(\alpha,L)$ is isomorphic, as a ${\mathbb Z}_2$-graded algebra, to the
Clifford algebra generated by the basis vectors $u,z_1,\ldots,z_n$
and with associated bilinear form corresponding to the matrix
$L=-\frac{1}{4\alpha}(2\alpha B+\Gamma)$. 
Since $\alpha$ is invertible, Proposition \ref{wk} in this case
 yields the well-known fact that a generalized Clifford algebra is
 ${\mathbb Z}_2$-graded central simple if and only if the associated
 bilinear form is
 non-degenerate. 

Let us recall the decomposition of $BM(\f, E(n), R_0)$ and  the
embedding of the subgroup $Sym_n(\f)$ described in \cite{GioJuan3}.
The pull-back along the injection $j$ of $\f{\mathbb Z}_2$ into $E(n)$ yields a
surjective map $j^*\colon BM(\f, E(n), R_0)\to BW(\f)$. The map $j^*$
is split, the
splitting map being induced by the pull-back $p^*$ along the projection
$p\colon E(n)\to \f{\mathbb Z}_2$. The Kernel of $j^*$ is isomorphic
to $Sym_n(\f)$ and representatives of its elements can be
constructed as follows. To a symmetric matrix $L$ 
one may associate a special $2$-cocycle $\omega$
such that
$$
\begin{array}{cc}
\omega(c\otimes c)=1;&\omega(x_i\otimes x_j)=-l_{ij};\\
\omega(c\otimes x_i)=0;&\omega(x_i\otimes c)=0.
\end{array}
$$
The left regular action of $_\omega E(n)$ given by
$f_h(k)=\sum\omega(\h1\otimes\k1)\h2\k2$  induces an inner action of
$E(n)$ on $A^{\omega}=End(_\omega E(n))$ by: $h\cdot f=\sum f_{\h1}\circ
f\circ f^{-1}_{\h2}$. The module algebra $A^\omega$ is $(E(n),
R_0)$-Azumaya and it represents an element in the kernel of $j^*$. The 
subalgebra $Ind(A^\omega)$ generated by $U=f_c$ and $W_i=-f^{-1}_{x_i}$ for
$i=1,\ldots,n$ is a submodule algebra, and its relations are:
$$
\begin{array}{ccc}
U^2=1,& W_iW_j+W_jW_i=2l_{ij},& UW_j+W_jU=0.\\
\end{array}
$$  
The action on $A^\omega$ is given by:
$$c\cdot f=U fU^{-1};\quad x_j\cdot f=W_j(c\cdot f)-fW_j$$
for every $f\in A^\omega$.
One shows that $E(n)$ acts innerly on any representative $A'$ of the class
of $A^\omega$. If the action is realized by the map $g\colon E(n)\to
A'$ with $g(c)$ and $g^{-1}(x_j)$ skew-commuting, then
the matrix $L$ describing the relations among the $g^{-1}(x_j)$'s is an
invariant of the class represented by  $A^\omega$ and it 
uniquely determines the class in the kernel of $j^*$.   

\smallskip

For every quasitriangular Hopf algebra $H$ and $(H, R)$-Azumaya algebra $A$ 
let us denote by $[A]$ the class in $BM(\f, H, R)$ represented by
$A$.  Let us define,
for any nonzero $t\in \f$, the algebra $C(t)$ generated by $x$ with
relation $x^2=t$, $c$-action $c\cdot x=-x$ and trivial action of the
$x_j$'s. Being the representative of an element of $BW(\f)$, 
the $E(n)$-module algebra $C(t)$ is $(E(n), R_0)$-Azumaya. It is well-known
that $\overline{C(t)}=C(-t)$. For evey symmetric $n\times n$ matrix
$L$ we denote by $A(1,L)$ the
$E(n)$-module algebra, representing $p^*\circ j^*([C(1,L)])$, 
isomorphic to $C(1,L)$ as a $\f{\mathbb Z}_2$-module
algebra and with trivial action of the $x_j$'s.

\begin{proposition}Let $\alpha$ be a nonzero element in $\f$ and let
$L$ be an $n\times n$ invertible, symmetric matrix with entries in
$\f$.  With notation as above
$$[C(\alpha,L)]=[C(-1)\#C(\alpha)\# A(1,L)][A^\omega]$$ with
$\omega(x_i\otimes x_j)=((4L)^{-1})_{ij}$.  
\end{proposition} 
\pf Let us first assume that $\alpha$ is a
square, so that $[C(\alpha, L)]=[C(1, L)]$.
The class 
$\left[C(1, L)\right]\left[\overline{A(1,L)}\right]=\left[C(1, L)\#\overline{A(1,L)}\right]$ belongs to the
Kernel of $j^*$. We compute its matrix invariant.

Let us denote the generators of $\overline{A(1,L)}$ by $U$ and
$Z_1,\,\ldots,\,Z_n$. The relations in
$\overline{A(1,L)}$ are:

$$
\begin{array}{cc}  
UZ_i+Z_iU=0;&Z_iZ_j+Z_jZ_i=-2l_{ij};\\
U^2=-1;&Z_i^2=-l_{ii}
\end{array}
 $$ and the action is determined by 
$$
\begin{array}{cc}
c\lu U=-U;& 
c\lu Z_i=-Z_i;\\
x_i\lu U=0;&x_j\lu
Z_i=0.
\end{array}
$$ 
In the product $C(1, L)\#\overline{A(1,L)}$ the elements $U, Z_j$ for
$j=1,\,\ldots,\,n$ 
skew-commute with $u$ and $z_j$ for $j=1,\,\ldots n$. Let us introduce the elements:
$$w_j=\sum ((-2L)^{-1})_{jk}z_k\in C(1,L).$$ It is not hard to see that 
$x_j\lu b=w_j(c\lu b)-bw_j$ for $b=u,w_j$ for $j=1,\ldots,n$. Since
this equality extends to products, it holds for every $b\in C(1,L)$. Besides, 
for every $a\#b\in C(1,L)\#\overline{A(1,L)}$ we have:
$$
\begin{array}{rl}
x_j\lu(a\#b)&=(x_j\lu a)\#(c\lu b)+a\#( x_j\lu b)\\
&=(x_j\lu a)\#(c\lu b)\\
&=w_j(c\lu a)\#(c\lu b)-aw_j\#c\lu b\\
&=(w_j\#1)(c\lu(a\# b))-(a\# b)(w_j\#1).
\end{array}
$$
Since $[C(1,L)\#\overline{A(1,L)}]\in Ker(j^*)$ 
the action of $c$ on this product is strongly inner, i.e.,
there exists an element $Y\in C(1,L)\#\overline{A(1,L)}$ with $Y^2=1$
and such that $c\lu(a\# b)=Y(a\# b)Y^{-1}$. In particular,
$Y(w_j\#1)+(w_j\#1)Y=0$. Therefore, the relations among the
$(w_j\#1)$'s for $j=1,\,\ldots, n$ will give the saught invariant matrix.    
This is easily computed and we have:
$$(w_i\#1)(w_j\#1)+(w_j\#1)(w_i\#1)=2((4L)^{-1})_{ij}.$$
Hence, 
\begin{equation}\label{deco}[C(t^2,L)]=[A(1,L)][A^\omega]\end{equation}
with $\omega(x_i\otimes x_j)=((4L)^{-1})_{ij}$. 

Let us assume now that $\alpha$ is not a square. As $E(n)$-module algebras
$C(1)\#C(\alpha, L)\cong C(\alpha)\#C(1,L)$. 
By formula (\ref{deco}) we have the statement.\hfill$\Box$

\smallbreak

The Hopf algebra corresponding to $n=1$ is just Sweedler's Hopf
algebra $H_4$. Proposition \ref{wk}, together with
self-duality of $H_4$, states that the module algebras
$A\langle\frac{\alpha,\,\beta,\,\gamma}{\f}\rangle$
 in \cite{yinhuo} are $(H_4,\,R_t)$-Azumaya if and only
if $2\alpha(t-2\beta)+\gamma^2\neq0$. If $\gamma=0$, we recover
 \cite[Proposition 3.1]{yinhuo2}. When $t=0$, we
recover the result in \cite{yinhuo} that the algebra
$A\langle\frac{\alpha,\,\beta,\,\gamma}{\f}\rangle$
is $(H_4,\,R_0)$-Azumaya if and only
if $-4\alpha\beta+\gamma^2\neq0$. Up to a slight change in notation,
Proposition \ref{deco} provides a bridge between the construction of
the map $(\f,+)\to BM(\f, H_4, R_0)$ in \cite{yinhuo} and the
construction of the map $Sym_n(\f)\to BM(\f, E(n),R_0)$ in
\cite{GioJuan3} for $n=1$.

\section*{Acknowledgements}

This research was partially supported by Progetto Giovani Ricercatori number
CPDG031245 of the University of Padova. The author wishes to thank Professor
Juan Cuadra of the University of Almer\'{\i}a for suggesting the problem
leading to Theorem \ref{analogo}.

\end{document}